\documentclass[12pt]{amsart}
\usepackage{amscd,amssymb}
\usepackage[arrow,matrix]{xy}
\usepackage[plainpages,backref,urlcolor=blue]{hyperref}

\theoremstyle{plain}
\newtheorem{thm}[subsection]{Theorem}
\newtheorem{lem}[subsection]{Lemma}
\newtheorem{prop}[subsection]{Proposition}

\newtheorem{conj}[subsection]{Conjecture}

\theoremstyle{definition}
\newtheorem{rk}[subsection]{Remark}

\newtheorem{ex}[subsection]{Example}
\newtheorem{question}[subsection]{Question}

\numberwithin{equation}{section}
\newcommand{\R}{\mathbb{R}}
\newcommand{\C}{\mathbb{C}}

\newcommand{\PP}{\mathbb{P}}

\newcommand{\N}{\mathbb{N}}

\begin{document} 

\title [Log-concavity of Milnor algebras for projective hypersurfaces]
{Log-concavity of Milnor algebras for projective hypersurfaces } 

\author[Gabriel Sticlaru]{Gabriel Sticlaru}
\address{Faculty of Mathematics and Informatics,
Ovidius University}
\email{gabrielsticlaru@yahoo.com } 



\begin{abstract}
We investigate the 0-th local cohomology of the Jacobian ring of a homogeneous polynomial defining a projective hypersurface whose singular locus is a 0-dimensional complete intersection.

\end{abstract}

\maketitle
 
\section{Introduction and main result}

A sequence $a_0,a_1, ..., a_m$ of  real numbers is said to be {\it log-concave} (resp. {\it strictly log-concave})
if it verifies $a_k^2 \geq a_{k-1}a_{k+1}$ (resp.    $a_k^2 > a_{k-1}a_{k+1}$)  for $k=1,2,...,m-1$. An infinite sequence $a_k$, $k \in \N$ is (strictly) log-concave if any truncation of it is (strictly) log-concave.
 Such sequences play an important role in Combinatorics and Algebraic Geometry, see for instance the recent paper
\cite{H}. Recall that a sequence $a_0,...,a_m$ of real numbers is said to be {\it unimodal} if there is an integer $i$ between $0$ and $m$ such that
$$a_0 \leq a_1\leq ...\leq a_{i-1} \leq a_i \geq a_{i+1} \geq ...\geq a_m.$$
A nonnegative log-concave sequence with no internal zeros (i.e. a sequence for which the indices of the nonzero elements form a set of consecutive integers), is known to be unimodal, see \cite{H} and the references there.

A polynomial $P(t) \in \R[t]$, 
$$P(t) = \sum_{j=0}^{j=m}a_jt^j,$$
 with coefficients $a_j$ for $0\leq j \leq m$, is said to be log-concave (resp. unimodal) if the sequence of its coefficients is log-concave (resp. unimodal).

Let $S=\C[x_0,...,x_n]$ be the graded ring of polynomials in $x_0,,...,x_n$ with complex coefficients and denote by $S_d$ the vector space of homogeneous polynomials in $S$ of degree $d$. 
For any polynomial $f \in S_r$ we define the {\it Jacobian ideal} $J_f \subset S$ as the ideal spanned by the partial derivatives $f_0,...,f_n$ of $f$ with respect to $x_0,...,x_n$. For $n=2$ we use $x,y,z$ instead of
$x_0, x_1, x_2$ and $f_x, f_y, f_z$  instead of $f_0, f_1, f_2$.

The Hilbert-Poincar\'{e} series of a graded $S$-module $M$ of finite type is defined by 
\begin{equation} 
\label{eq2}
HP(M)(t)= \sum_{k\geq 0} \dim M_k\cdot t^k .
\end{equation} 
For any polynomial $f \in S_d$ we define the hypersurface $V(f)$ given by $f=0$ in $\PP^n$ and the corresponding graded {\it Milnor} (or {\it Jacobian}) {\it algebra} by
\begin{equation} 
\label{eq1}
M(f)=S/J_f.
\end{equation}
Smooth hypersurfaces $V(f)$ of degree $d$ have all the same  Hilbert-Poincar\'{e} series, namely
\begin{equation} \label{eq1.5}
HP(M(f))(t)=\frac{(1-t^{d-1})^{n+1}}{(1-t)^{n+1}}=(1+t+t^2+\ldots + t^{d-2} )^{n+1} =\sum_{k=0}^{k=T}a_{k}t^{k},
\end{equation} 
where $T=(n+1)(d-2)$ and the corresponding sequence of coefficients $a_0,...,a_T$ is clearly
log-concave. As soon as the hypersurface $V(f)$ acquires some singularities, the series $HP(M(f))$ is an infinite sum, and the sequence of coefficients $a_0,...,a_T$ may no longer be
log-concave, see Proposition 3.1 in \cite{St}.

For a hypersurface $V(f):f=0$ with isolated singularities we introduce (see \cite{DS2} for other invariants) the {\it coincidence threshold} defined as:
$$ct(V(f))=\max \{q~~:~~\dim M(f)_k=\dim M(f_s)_k \text{ for all } k \leq q\},$$
with $f_s$  a homogeneous polynomial in $S$ of degree $d=\deg f$ such that $V(f_s):f_s=0$ is a smooth hypersurface in $\PP^n$.

\bigskip

E. Sernesi in the recent preprint \cite{S1} suggests that many properties of the Milnor algebra $M(f)$ of a smooth hypersurface are replced by properties of the local cohomology module
\begin{equation} 
\label{eq2}
N(f)=H^0_{\bf m}(M(f))= \widehat J_f/J_f,
\end{equation}
where ${\bf m}$ is the maximal ideal $(x_0,x_1,...,x_n)$ in $S$ and $\widehat J_f$ is the saturation of the Jacobian ideal $J_f$, see \cite{D1}. This module was considered also in  \cite{DS},  \cite{SW} and it is closely related to the syzygies of the Jacobian ideal $J_f$ considered in  \cite{DS1},  \cite{DS2},  \cite{DS3}. 

In this note we consider the case when the projective hypersurface $V(f)$ has only isolated singularities, which implies that the corresponding series $HP(N(f))(t)$ is a polynomial.
We propose the following conjecture, which is a stronger form of Dimca's conjecture in Remark 6
in  \cite{D1}.
\begin{conj} \label{conj1}
 If the projective hypersurface $V(f)$ has only isolated singularities, then the corresponding series $HP(N(f))(t)$ is a log-concave polynomial with no internal zeros. In particular, $HP(N(f))(t)$ is a unimodal polynomial.
\end{conj} 
The difficulty to approach this conjecture lies in the fact that there is too little information on the coefficients $n_k=\dim N(f)_k$ in general. There is however a case when one has complete information on these coefficients. Indeed, when the saturated Jacobian ideal $\widehat J_f$ is a complete intersection ideal of multidegree $(d_1,...,d_n)$, then it is known that
\begin{equation} 
\label{eq3}
HP(M(f))(t)=\frac{(1-t^{d-1})^{n+1}+t^{(n+1)(d-1)-\sum d_i} (1-t^{d_1}) \cdots (1-t^{d_n})}{(1-t)^{n+1}}.
\end{equation}
Also, we have Tjurina number $\tau(V(f))=d_1\cdots d_n$ and $ct(V(f))=T-\sum d_i +n$, see Proposition 4. in  \cite{D1}. Then
$HP(N(f))(t)=HP(S/J_f)(t)-HP(S/\widehat J_f)(t)$ and hence
\begin{equation} 
\label{eq4}
HP(N(f))(t)=\frac{(1-t^{d-1})^{n+1}}{(1-t)^{n+1}}-\frac{1-t^{(n+1)(d-1)-\sum d_i}} {1-t} \cdot
\frac{1-t^{ d_1}} {1-t} \cdots \frac{1-t^{ d_n}} {1-t}.
\end{equation}
In spite of this very explicit formula, we are not able to prove the conjecture in this case for arbitrary $n$. 

The main result of this note is the following.

\begin{thm} \label{thm1}
The Conjecture \ref{conj1} is true for polynomials $f$ such that the saturated Jacobian ideal $\widehat J_f$ is a complete intersection ideal and $n=2$. 
\end{thm}
There is also the following related question, which is a stronger form of Dimca and Saito's Question 1 in \cite{DS}.

\begin{question} \label{q1}
 If the projective hypersurface $V(f)$ has only isolated singularities, then the difference $HP(M(f_s))(t)-HP(N(f))(t)$ is log-concave and unimodal polynomial, where $f_s\in S$ is any homogeneous polynomial of degree $d$ in such that $V(f_s)$ is smooth.
\end{question} 
The answer is far from being known in general, but the formulas in \eqref{eq1.5} and \eqref{eq4} imply the following result for any $n$.
\begin{prop} \label{prop1}
The answer to Question \ref{q1} is positive for polynomials $f$ such that the saturated Jacobian ideal $\widehat J_f$ is a complete intersection ideal. 
\end{prop}
It is a classic result (see for instance \cite{B1} for a proof) that the product of log-concave polynomials is log-concave.

Unlike Proposition \ref{prop1}, whose proof is straightforward, the proof of Theorem \ref{thm1} is rather long and tedious. In the next section we obtain formulas for the coefficients of polynomial $HP(N(f))(t)$ via a key result due to  D. Popescu and M. Vl\u adoiu in \cite [Lemma 2.9]{PV}.
The proof is divided in eight subcases, and we show by simple examples that all these subcases really arise for some singular plane curves.
In the final section we check the log-concavity in all these cases.

\begin{rk} \label{rkConca}

Recently Professor Aldo Conca has communicated us the following two examples with non log-concave and/or non-unimodal Hilbert-Poincar\'e series.
In these examples the saturated
Jacobian ideal is not a complete intersection ideal.

(i) This is an example of non log-concave projective curve in $\PP^2$, namely
$$C: f=x^{9}y+y^{10}+x^{3}y^{5}z^{2}=0.$$
Here $$HP(N(f))(t)=t^{8}+2t^{9}+3t^{10}+5t^{11}+6t^{12}+5t^{13}+3t^{14}+2t^{15}+t^{16}$$
is not a log-concave polynomial because the sequences of coefficients $ 2,3,5$ is not log-concave.

(ii) This is an example of non unimodal projective surface in $\PP^3$, namely
$$D: f=x^{5}+y^{5}+yz(x^{3}+z^{2}w)=0.$$
Here $HP(N(f))(t)=t^{5}+t^{7}$ is  not unimodal because the sequences of coefficients $1,0,1$ has internal zeroes.

In conclusion, we must add some  restrictions on $f$ such that the  Conjecture \ref{conj1} may be true. This is exactly what we have done in Theorem \ref{thm1} above.

\end{rk}

The author would like to thank Alexandru Dimca for suggesting this interesting problem and to Aldo Conca for the surprizing examples above.
\section{Determining the coefficients}
From now on we suppose $n=2$ and set $a=d_1$ and $b=d_2$ to simplify the notation. With this notation we have
\begin{equation} 
\label{eq5}
HP(N(f))(t)=\frac{(1-t^{d-1})^{3}}{(1-t)^{3}}-\frac{1-t^{3d-3-a-b}} {1-t} \cdot
\frac{1-t^{ a}} {1-t} \cdot \frac{1-t^{ b}} {1-t},
\end{equation}
or, equivalently,
\begin{equation} 
\label{eq6}
HP(N(f))(t)=P_{d-1}^3(t)-P_a(t) P_b(t)P_{3d-3-a-b}(t),
\end{equation}
with $P_m(t)=(1-t^m)/(1-t)$. We can and will assume in the sequel that 
$$1 \leq a \leq b <d.$$
One clearly has $HP(N(f))(t)=t^aQ(t)$, where $Q(t)$ is a symmetric polynomial of degree $T-2a=3d-6-2a$, unless $a=b=d-1$ when $Q(t)=0$.

\bigskip

\noindent The coefficients of a general product
$$A(t)=(1+t+t^2+\ldots +t^{a-1}) (1+t+t^2+\ldots +t^{b-1})(1+t+t^2+\ldots +t^{c-1})$$
are described by D. Popescu and M. Vl\u adoiu in \cite [Lemma 2.9]{PV} and this results plays a key role in this section.

\bigskip

\noindent If we set $a=b=c=d-1$, this yields the coefficients of the polynomial
 $A(t)=(1+t+t^2+\ldots +t^{d-2})^3$ in the following form
\[ a_k = \begin{cases}
	    \binom{k+2}{2}, 0 \leq k \leq d-2 \\
	   \binom{d}{2}+\sum _{i=1} ^{k-d+2} (d-1 -2i),  k \geq d-1 . 
	 \end{cases}
\] 
Equivalently, one has
\[ a_k = \begin{cases}
	    \binom{k+2}{2}, 0 \leq k \leq d-2 \\
	   \binom{k+2}{2}-3\binom{k+3-d}{2}, d-1 \leq k \leq T/2 \\
      a_{T-k}, k>T/2
	 \end{cases}
\] 
Using \cite [Lemma 2.9]{PV}  we determine now the coefficients $q_k$ of the polynomial  $Q(t)=q_0+q_1t+q_2t^2+\ldots  $   for $0 \leq k \leq T/2-a$, which is enough by symmetry. 
The number of cases being quite large, we leave the proofs for the reader.

\subsection{{Case 1: $2a+2b+1<3d, a \leq b < d$}}

\vskip1cm

\center{\it{Case 1.1: $2a+2b+1<3d, a=b<d$}}

\bigskip
\begin{itemize}
\item  subcase (1.1.1):      $2a+2b+1<3d, 2 \leq a=b<d, a \leq d-2 \leq 2a-2$.
\end{itemize} 
\medskip 
\[ q_k = \begin{cases}
	    2\binom{k+2}{2}, 0 \leq k \leq d-a-2 \\
	   2\binom{k+2}{2}-3\binom{k+a+3-d}{2}, d-a-1 \leq k \leq a-2 \\
	   \binom{k+a+2}{2}-3\binom{k+a+3-d}{2} - a^2,  a-1 \leq k \leq T/2-a
	 \end{cases}
\]

\medskip
\begin{itemize}
\item subcase (1.1.2):   $2a+2b+1<3d, 1 \leq a=b<d, d-2 \geq 2a-1$.
\end{itemize} 
\medskip 
\[ q_k = \begin{cases}
	    2\binom{k+2}{2}, 0 \leq k \leq a-2 \\
	   \binom{k+a+2}{2}-a^2, a-1 \leq k \leq d-a-2 \\
	   \binom{k+a+2}{2}-3\binom{k+a+3-d}{2} - a^2,  d-a-1 \leq k \leq T/2-a
	 \end{cases}
\]

\center{\it{Case 1.2: $2a+2b+1<3d, a<b<d$}}

\medskip
\begin{itemize}
\item subcase (1.2.1): $2a+2b+1<3d, 1 \leq a<b<d, a \leq d-2 \leq b-1$, i.e. $d=b+1$.
\end{itemize} 
\medskip 
\[ q_k = \begin{cases}
	    \binom{k+2}{2}, 0 \leq k \leq b-a-1 \\
	   \binom{k+2}{2}-2\binom{k+2-(b-a)}{2}, b-a \leq k \leq b-2 \\
	   \binom{k+a+2}{2}-3\binom{k+2-(b-a)}{2} - ab,  b-1 \leq k \leq T/2-a
	 \end{cases}
\]

\medskip
\begin{itemize}
\item subcase (1.2.2): $2a+2b+1<3d, 2 \leq a<b<d, b \leq d-2 \leq a+b-2$.
\end{itemize} 
\medskip 
\[ q_k = \begin{cases}
	    \binom{k+2}{2}, 0 \leq k \leq b-a-1 \\
	   \binom{k+2}{2}+\binom{k+2-(b-a)}{2}, b-a \leq k \leq d-a-2 \\
	   \binom{k+2}{2}+\binom{k+2-(b-a)}{2}-3\binom{k+a+3-d}{2},  d-a-1 \leq k \leq b-2 \\
   \binom{k+a+2}{2}-3\binom{k+a+3-d}{2} - ab,  b-1 \leq k \leq T/2-a
	\end{cases}
\]

\medskip
\begin{itemize}
\item  subcase (1.2.3): $2a+2b+1<3d, 1 \leq a<b<d, d-2 \geq a+b-1$.
\end{itemize} 
\medskip 
\[ q_k = \begin{cases}
	    \binom{k+2}{2}, 0 \leq k \leq b-a-1 \\
	   \binom{k+2}{2}+\binom{k+2-(b-a)}{2}, b-a \leq k \leq b-2 \\
	   \binom{k+a+2}{2} - ab,  b-1 \leq k \leq d-a-2 \\
   \binom{k+a+2}{2}-3\binom{k+a+3-d}{2} - ab,  d-a-1 \leq k \leq T/2-a	
	\end{cases}
\]

\subsection{{Case 2: $2a+2b+1 \geq 3d, a \leq b < d$}}

\center{\it{Case 2.1: $2a+2b+1 \geq 3d, a=b<d$}}

\medskip
\begin{itemize}
\item subcase (2.1.1): $2a+2b+1 \geq 3d, 2 \leq a=b<d, a \leq d-2$.
\end{itemize} 

\medskip 
\[ q_k = \begin{cases}
	    2\binom{k+2}{2}, 0 \leq k \leq d-a-2 \\
	   2\binom{k+2}{2}-3\binom{k+a+3-d}{2}, d-a-1 \leq k \leq 3d-3a-6 \\
	   3(d-a-1)^2,  3d-3a-5 \leq k \leq T/2-a
	\end{cases}
\]

\center{\it{Case 2.2: $2a+2b+1 \geq 3d, a<b<d$}}
\medskip
\begin{itemize}
\item subcase (2.2.1): $2a+2b+1 \geq 3d, 2 \leq a<b<d, a \leq d-2 \leq b-1$, i.e. $d=b+1$.
\end{itemize} 
\medskip 
\[ q_k = \begin{cases}
	    \binom{k+2}{2}, 0 \leq k \leq b-a-1 \\
	   \binom{k+2}{2}-2\binom{k+2-(b-a)}{2}, b-a \leq k \leq 2b-2a-3 \\
	   (b-a)^2, 2b-2a-2 \leq k \leq T/2-a 
	\end{cases}
\]

\medskip
\begin{itemize}
\item subcase (2.2.2): $2a+2b+1 \geq 3d, 2 \leq a<b<d, b \leq d-2$.
\end{itemize} 
\medskip 
\[ q_k = \begin{cases}
	    \binom{k+2}{2}, 0 \leq k \leq b-a-1 \\
	   \binom{k+2}{2}+\binom{k+2-(b-a)}{2}, b-a \leq k \leq d-a-2 \\
     \binom{k+2}{2}+\binom{k+2-(b-a)}{2}-3\binom{k+a+3-d}{2}, d-a-1 \leq k \leq 3d-2a-b-6 \\  
     (b-a)^2+3(d-a-1)(d-b-1), 3d-2a-b-5 \leq k \leq T/2-a 
	\end{cases}
\] 

\begin{rk} \label{rk1}
 It is likely that there are geometric restrictions on the possible values of $a$ and $b$ in relation to the fact that the corresponding polynomials  span the saturated Jacobian ideal $\widehat J_f$ of a plane curve $V(f)$.
For instance, it follows from Corollary 5 in  \cite{D1} that $a+b=T+2-ct(V(f))$, and hence the condition in Case 2, i.e. $2a+2b+1 \geq 3d$ is equivalent to $ct(V(f))<T/2$.

However, we present  in the following table,  for each subcase  above an example of a curve for which the saturated Jacobian ideal $\widehat J_f$ is a complete intersection satisfying the inequalities defining the subcase. 
 
\noindent
\begin{center}
 \begin{tabular} {|l|c|l|}
  \hline
  	          &   Configuration     &                   \\                      
  Subcase     &      $(a,b,d)$      &   polynomial   \  \  $ f$  \\ 
  \hline
 
  $ 1.1.1 $ &  $ (2, 2, 4)   $  & $ x^2y^2+xz^3+yz^3                       $ \\
  $ 1.1.2 $ &  $ (1, 1, 4)   $  & $ xyz^2+x^4+y^4                          $ \\
  $ 1.2.1 $ &  $ (4, 7, 8)   $  & $ x^3y^5+y^8+z^8                         $ \\
  $ 1.2.2 $ &  $ (3, 4, 6)   $  & $ (x^2+y^2)^3+(y^3+z^3)^2  							 $ \\
  $ 1.2.3 $ &  $ (2, 3, 6)   $  & $ x^2(x+z)^2(x-z)^2-y^2(y-z)^2(y^2+2z^2) $ \\
  $ 2.1.1 $ &  $ (5, 5, 7)   $  & $ z(x^6+y^6)+2y^7                        $ \\
  $ 2.2.1$ &   $ (5, 6, 7)   $  & $ x^4y^3+z^7                             $ \\
  $ 2.2.2 $ &  $ (8, 10, 12) $  & $ (x^4+y^4)^3+(y^2+z^2)^6                $ \\

  \hline
  
\end{tabular} 
\end{center}

\end{rk}

 \begin{ex}
\label{ex1}

Let $V(f): f=f_1^2+f_2^2=0$ where $f_1=x^2+y^2, f_2=y^2+z^2$. This is a nodal curve with degree $d=deg(f)=4$ and four singularities of type $A_1,$ i.e. nodes: $(1:i:1), (1:-i:1), (1:i:-1), (1:-i:-1).$ 
The Jacobian ideal of $f$ is $J_f=(x^3+xy^2, x^2y+2y^3+yz^2, y^2z+z^3)$ and the saturated Jacobian ideal is 
$\widehat J_f=(y^2+z^2, x^2-z^2)=(f_1, f_2)$ which implies that $\widehat J_f$ is a complete intersection ideal of multidegree $(a=2,b=2).$ 

$HP(M(f_s))(t)=\frac{(1-t^3)^3}{(1-t)^3}=(1+t+t^2)^3 =1+3t+6t^2+7t^3+6t^4+3t^5+t^6.$

$HP(M(f))(t)= \frac{(1-t^3)^3+t^5 (1-t^2) \cdot (1-t^2)}{(1-t)^3} =1+3t+6t^2+7t^3+6t^4+4(t^5+ \cdots.$

$HP(N(f))(t)=\frac{(1-t^3)^3}{(1-t)^3}-\frac{1-t^5} {1-t} \cdot
\frac{1-t^2} {1-t} \cdot \frac{1-t^2} {1-t} = 2t^2+3t^3+2t^4 $, 
a log-concave and unimodal polynomial. 

$HP(M(f_s))(t)-HP(N(f))(t)= \frac{1-t^5} {1-t} \cdot
\frac{1-t^2} {1-t} \cdot \frac{1-t^2} {1-t} = 1+3t+4t^2+4t^3+4t^4+3t^5+t^6$, symmetric, log-concave and unimodal polynomial.

This is the subcase $(1.1.1):$ $2a+2b+1<3d, 2 \leq a=b<d, a \leq d-2 \leq 2a-2$.
More generally, if $f=f_1^2+f_2^2$ with $f_1=x^m+y^m, f_2=y^m+z^m$, the saturated Jacobian ideal is 
$\widehat J_f=(y^m+z^m, x^m-z^m)=(f_1, f_2)$ which implies that $\widehat J_f$ is a complete intersection ideal of multidegree $(a=m,b=m).$ 

\end{ex}

\begin{ex}
\label{ex2}

Let $V(f): f=(x^2+y^2)^3+(y^3+z^3)^2=0$ of type $6A_2$, with six cusp as singularities.  The saturated Jacobian ideal is $\widehat J_f=(y^3+z^3, x^4+2x^2y^2-yz^3),$ with two generators which implies that $\widehat J_f$ is a complete intersection ideal of multidegree $(a=3,b=4).$ 

This is the subcase $(1.2.2):$ $2a+2b+1<3d, 2 \leq a<b<d, b \leq d-2 \leq a+b-2$.

If the saturated ideal is difficult to compute, we can try to decide if the saturated Jacobian ideal is a complete intersection based on Corollary $5$ from \cite{D1}: when $n=2$, the couple $(a,b)$ is determined, when it exists, by the couple $(\tau(V(f)),ct(V(f)))$.

We have Tjurina number $\tau(V(f))=a\cdot b=12$ and $ct(V(f))=T-(a+b) +n=7$ so $a=3, b=4$.

Hilbert-Poincar\'{e} series is 
$HP(M(f))(t)=1+3t+6t^2+10t^3+15t^4+18t^5+19t^6+18t^7+16t^8+13t^9+12(t^{10}+ \cdots .$

On the other hand, if we assume the complete intersection case,  with the formulas in  \eqref{eq3} we get the same series: $\frac{(1-t^5)^3+t^8 (1-t^3) \cdot (1-t^4)}{(1-t)^3}$. So, the conclusion of such a computation is that the saturated Jacobian ideal might be a complete intersection

\end{ex}

\begin{ex}
\label{ex3}

 Consider now the irreducible curve $V(f):x^3z^4+xy^5z+x^7+y^7=0$, which has one {\it non weighted homogeneous singularity} located at $(0:0:1)$ with Milnor number $\mu(V(f))= 12$, Tjurina number $\tau(V(f))=11$ and $ct(V(f))=11$.
The saturated Jacobian ideal is $\widehat J_f=(x^3, x^2y^2, y^5+3x^2z^3, 5xy^4-21x^2yz^2)$ with four polynomials as generators, which implies that $\widehat J_f$ is not a complete intersection ideal. For this case, the system: $a\cdot b=11$ and $a+b=6$ has no solutions.

The computation using CoCoA \cite{Co} or Singular \cite{Sing} yields in this case \\
$HP(M(f))(t)=1+3t+6t^2+10t^3+15t^4+21t^5+25t^6+27t^7+27t^8+25t^9+21t^{10}+15 t^{11}+12t^{12}+11(t^{13} \cdots $  \\ 
$HP(N(f))(t)=t^3+4t^4+10t^5+14t^6+16t^7+16t^8+14t^9+10t^{10}+4t^{11}+t^{12}.$ \\
$HP(M(f_s))(t)-HP(N(f))(t)=
1+3t+6t^2+9t^3+11t^4+11t^5+11t^6+11t^7+11t^8+11t^9+11t^{10}+11t^{11}+9t^{12}+6t^{13}+3t^{14}+t^{15}.$

This is not complete intersection case but the Conjecture \ref{conj1} and the Question \ref{q1} are true. 

\end{ex}

\section{Checking the log-concavity}

To prove the log-concavity of the polynomial $Q(t)$ we use the following lemmas.

\begin{lem} \label{lem1} 
The following sequences are strictly log-concave. 

\begin{enumerate}
\item	A lstrictly og-concave sequence multiplied by a strictly positive real number.

\item	$a_k=\binom{m+k}{2}$ for $k \in \N$.
\item	$a_k=\binom{m+k}{2}-n$ for $k \in \N, n>0$.
\item	$a_k=x\binom{m+k}{2}-y\binom{n+k}{2}$, for $x, y > 0$ and $x/y<24$. 
\item	$a_k=\binom{m+k}{2}+\binom{n+k}{2}$, for $m, n \geq 2, k > 0$. 

\end{enumerate}

\end{lem} 
\noindent Note that the sum of two log-concave sequences (resp. the difference between a 
 log-concave sequence and a constant sequence) is not in general log-concave.

\proof

All the claims are easy to check, we give the details only for the case (4).
Set $ a=m+k $ and $ b=n+k $ and note that 
$$a_k^2-a_{k-1}a_{k+1}=1/2(a^2x^2+a^2xy-4abxy+b^2xy+b^2y^2-ax^2+axy+bxy-by^2).$$
This quadratic function in $a$ has the discriminant
$$\Delta= -x(x-y)^2(4b^2y-4by-x)$$
which is negative for $x/y<4b(b-1)$. Since $b-1 \geq 2$, it follows that for $x/y<24$ we have log-concavity.

\endproof

\begin{lem} \label{lem2} 
The following sequences truncated at $k\leq T/2-a$ are strictly log-concave.

\begin{enumerate}
\item	The sequence in the subcase (1.1.1) given by
$q_k= \binom{k+a+2}{2}-3\binom{k+a+3-d}{2} - a^2$, for $ a-1 \leq k $ , $a+2 \leq d \leq 2a$.
	  
\item	The sequence in the subcase (1.1.2) given by 
$q_k= \binom{k+a+2}{2}-3\binom{k+a+3-d}{2} - a^2$,  for $d-a-1 \leq k$ , $d \geq 2a+1$.

\item	The sequence in the subcase (1.2.1) given by  
 $q_k=  \binom{k+a+2}{2}-3\binom{k+2-(b-a)}{2} - ab,$ for $  b-1 \leq k$ , $b \geq 2a-1$.

\item	The sequence in the subcases (1.2.2) and (2.2.2) given by  
$q_k=  \binom{k+2}{2}+\binom{k+2-(b-a)}{2}-3\binom{k+a+3-d}{2},$ for $  d-a-1 \leq k$,   $b+2 \leq d $.

\item	The sequence in the subcase (1.2.2) given by  
$q_k= \binom{k+a+2}{2}-3\binom{k+a+3-d}{2} - ab$,  for $b-1 \leq k$,  $b+2 \leq d \leq a+b$.

\item	The sequence in the subcase (1.2.3) given by  
$q_k= \binom{k+a+2}{2}-3\binom{k+a+3-d}{2} - ab$,  for $d-a-1 \leq k$,  $ d \geq a+b+1$.

\end{enumerate}

\end{lem} 

\proof

Again all the claims are easy to check, we give the details only for the case (4) as a sample.
We set $k=d-a-1+x, d=b+2+y, a=1+z$ and we get  
$$q_{k}^2-q_{k-1}q_{k+1}=(x^2-(4y+2z+3)x+(10y^2+10yz+3z^2+24y+13z+16))/2 >0 $$
since the discriminant with respect to $x$,
$\Delta =-(6y^2+6yz+2z^2+18y+10z+55/4)  $  is negative.

\endproof
The above two lemmas imply that all the subsequences given by the same formulas are log-concave. For instance, in the subcase (1.1.1), the three subsequences
 $$2\binom{k+2}{2}, 0 \leq k \leq d-a-2,$$
	   $$2\binom{k+2}{2}-3\binom{k+a+3-d}{2}, d-a-1 \leq k \leq a-2 $$ and
	   $$\binom{k+a+2}{2}-3\binom{k+a+3-d}{2} - a^2,  a-1 \leq k \leq T/2-a$$
are log-concave. It remains to check what happens at the nodes, i.e. at values of $k$ where the formulas change, i.e. in the above cases we have to consider the sign of the difference $q_{k}^2-q_{k-1}q_{k+1}$ for  $k=d-a-2, d-a-1,a-2, a-1$.

\begin{lem} \label{lem3} The inequality $q_{k}^2-q_{k-1}q_{k+1} > 0$ holds at any node
in the subcases listed in section 2.
\end{lem}
\proof
As a sample we treat only the four values of $k$ above associated with the subcase (1.1.1).
We set $d=(4a+1)/3+1+x,a=1+y$ with $x,y \le 0$.

\noindent
$q_{d-a-2}^2 - q_{d-a-3} q_{d-a-1} =  \left[ 2\binom{d-a}{2} \right]^2 - 2\binom{d-a-1}{2} \left[ 2\binom{d-a+1}{2}-3\binom{2}{2} \right]=(3x+y+2)(15x+5y+7)/9 >0$.

\noindent
$q_{d-a-1}^2 - q_{d-a-2} q_{d-a} =   \left[ 2\binom{d-a+1}{2}-3\binom{2}{2} \right]^2
-2\binom{d-a}{2} \left[ 2\binom{d-a+2}{2}-3\binom{3}{2} \right]=(9x^2+6xy+y^2+39x+13y+13)/9 >0.$

\noindent
$q_{a-2}^2 - q_{a-3} q_{a-1} =$

$
\left[ 2\binom{a}{2}-3\binom{2a-d+1}{2} \right]^2 -
\left[ 2\binom{a-1}{2}-3\binom{2a-d}{2} \right]
\left[ 2\binom{2a+1}{2}-3\binom{2a-d+2}{2}-a^2 \right] = (45x^2+12xy+2y^2+51x+8y+14)/6 >0$. 

\noindent
$q_{a-1}^2 - q_{a-2} q_{a} =$

$\left[ 2\binom{2a+1}{2}-3\binom{2a-d+2}{2}-a^2 \right]^2 -
\left[ 2\binom{a}{2}-3\binom{2a-d+1}{2} \right]
\left[ 2\binom{2a+2}{2}-3\binom{2a-d+3}{2}-a^2 \right]$ 

$=(18x^2+12xy+2y^2+15x+8y+5)/3 >0$.

The above three Lemmas complete the proof of Theorem \ref{thm1}. Indeed, the inequality
$q_{k}^2-q_{k-1}q_{k+1} > 0$ holds unless the subsequence $q_{k-1},q_k, q_{k+1}$ is  constant, i.e. contained in the constant sequences occuring in the subcases (2.1.1), (2.2.1) and (2.2.2). But such constant sequences involve only stricty positive numbers. As $q_k \geq 0$ for any $k$, as it comes from the Hilbert-Poincar\'e series of $N(f)$, it follows that the sequence $q_k$ has no internal zeros, and hence it is unimodal as explained in the Introduction.


\end{document}